\documentclass{article}
\usepackage{ucs}
\usepackage[greek,english,spanish]{babel}
\usepackage{amsmath,amsfonts,amssymb}
\usepackage{textcomp}
\usepackage{wasysym}
\usepackage{calc}
\usepackage{longtable}
\usepackage{lscape}
\usepackage{supertabular}
\usepackage[labelfont=bf]{caption}
\usepackage{gensymb}
\usepackage{MnSymbol} 
\usepackage{textcomp} 
\usepackage[normalem]{ulem}
\usepackage{graphics}
\usepackage{float}
\usepackage[titletoc]{appendix}
\usepackage{hyperref}
\usepackage{tabulary}
\usepackage{array}
\usepackage{hhline}
\usepackage{longtable}
\usepackage[square,numbers]{natbib}
\usepackage{lineno}
\usepackage{url}

\usepackage{setspace} 
\usepackage{hypernat}

\newtheorem{theorem}{Theorem}

\newtheorem{lemma}{Lemma}

\hypersetup{colorlinks=true, linkcolor=blue, filecolor=blue, urlcolor=blue}

\newcommand{\Var}{\mathrm{Var}}

\newcommand\textstyleEmphasis[1]{\textit}

\makeatletter

\newcommand\ps@Standard{%
\renewcommand\@oddhead{}%
\renewcommand\@evenhead{}%
\renewcommand\@oddfoot{\thepage}%
\renewcommand\@evenfoot{\@oddfoot}%
\setlength\paperwidth{8.5in}
\setlength\paperheight{11in}
\setlength\voffset{-1in}
\setlength\hoffset{-1in}
\setlength\topmargin{0.7874in}
\setlength\headheight{12pt}
\setlength\headsep{0cm}
\setlength\footskip{12pt+0.1965in}
\setlength\textheight{11in-0.7874in-0.7874in-0cm-12pt-0.1965in-12pt}
\setlength\oddsidemargin{0.7874in}
\setlength\textwidth{8.5in-0.7874in-0.7874in}
\renewcommand\thepage{\arabic{page}}
\setlength{\skip\footins}{0.0398in}
\renewcommand\footnoterule{\vspace*{-0.0071in}
\noindent\textcolor{black}{\rule{0.25\columnwidth}{0.0071in}}
\vspace*{0.0398in}}
}

\newcommand{\abs}[1]{\lvert#1\rvert}

\renewcommand\paragraph{\@startsection{paragraph}{4}{\z@}%
{-2.5ex\@plus -1ex \@minus -.25ex}%
{1.25ex \@plus .25ex}%
{\normalfont\normalsize\itshape}}
\makeatother
\title{Fractal analysis of Pi normality}

\author{Carlos Sevcik\footnote{{Phone: +58 212 5051399. Mobile: +58 412 9319162. E-mail; csevcik@ivic.gob.ve. Private e-mail: carlos.sevcik.s@gmail.com}}}

\begin{document}
	\selectlanguage{english}
	\maketitle
	
	Laboratory on Cellular Neuropharmacology, {\selectlanguage{spanish}Centro de Biof\'{\i}sica y Bioqu\'{\i}mica, Instituto Venezolano de Investigaciones Cient\'{\i}ficas (IVIC), Apartado 20632, Caracas 1020A, Venezuela.}

\selectlanguage{english}

\begin{abstract}
$\pi$, the ratio between a circumference of a circle and is diameter, is a transcendental (and hence irrational) number.  Fractal analysis is used here to show that $\pi$\textquoteright{s} digit sequence corresponds to a uniformly distributed random succession of independent decimal digits, and that these properties get clearer as the number of digits in the series grows towards infinity; $ 10^9 $ digits were tested in this work. This is consistent with the hypothesis that $ \pi $ is normal.
\end{abstract}

\section{Introduction.} \label{sec:Intro}

The ratio of the circumference of a circle to its diameter, $\pi$, is a transcendental numbe \citep{Lambert1772, Kempner1916}. It has been proposed that $\pi$ is a normal number, although a definitive proof of this is lacking. A formal definition of number normality is in  \citet[Chapter 7.5, page 299]{Sierpinski1988}, but informally, a number $ w $ is normal to base $b$ if every sequence of $k$ consecutive digits in the base-$b$ expansion $ w $ appears with limiting frequency $b^{-k}$.  In other words, if a constant is normal to base 10, its decimal expansion would exhibit a \textquotedblleft{}7\textquotedblright{} one-tenth of the time, the string \textquotedblleft{}37\textquotedblright{} one one-hundredth of the time, and so on \citep{Bailey2001}. This means that every integer $\{c_1, c_2,\ldots,c_i,\ldots,c_b\}$ is equally likely to appear in the digit sequence $\{w_i\} $ of length $N$ corresponding to $ w $ in base $ b $. Probabilistically this means that all $ c_i $ are distributed according to a uniform (some times called rectangular) probability density function (\textit{pdf}) $ U_{\mathbb{Z}}(0, b-1) $, \citep{Wilks1962} and that
\begin{equation}
\left( w_k  \upmodels  w_{j \neq k}\right)  \forall\left( w_{j \neq k} \in\ \mathbb{Z} \right)_{\substack{ j=1,\ldots, N \\
k=1,\ldots,N}}
\end{equation}\label{eq:IntSec}
where $ \upmodels $ indicates stochastic independence \citep{Spohn1980}, and is opposite to $ \nupmodels $ meaning stochastic dependence. Tests for normality of $ \pi $ have been of a statistical nature \citep{Bailey2012a}. This analysis is by necessity limited to relatively simple combinations, usually less than five digits; but the $\pi$ digits sequence has been recorded to 22,459,157,718,361 decimal and 18,651,926,753,033 hexadecimal digits \citep{Trueb2016a, Trueb2016b}. Thus complex independent structures could remain hidden within the $ \pi $ digit series. Fractal analysis, on the other hand, treats an object as a whole, not just a fraction of it. Even if complexity is a difficult concept to define, fractal analysis provides insight on complex system in a unique manner which statistics or Euclidean geometry concepts cannot reach \citep{Mandelbrot1983, Mandelbrot1986, Barnsley1993, Hastings1993,  DSuze2105a}.

The uniform or rectangular pdf $ U[c,d] $, has no mode; $\forall \{c \leqq u  \leqq d\} \in \mathbb{R} $ have the same probability of occurring: $ (d-c)^{-1} $ \citep{Wilks1962}. All $ u $ are equally likely, but the next random variable produced by a random process obeying $ U[c,d] $ cannot be predicted. All the moments of the pdf are known but the information about the next event is nil. Classical \citet{Boltzmann1866, Boltzmann1964} entropy is $ S = k_B \, \ln(W) $, where $ k_B $ is the Boltzmann constant and $ W $ is the number of system states. If the set $ \{u\}_{i=1,2,\ldots, \infty} $ is a sequence of all randomly sampled $ u \in [c,d] $ any permutation of the $ \{u\}_i $ set is equally likely, the superset $  \{u\}_i \supseteq \{\Omega\} $ contains  $ W=\infty $ $ \{u\}_i $ and $ S_{\{\Omega\}}= \infty $. In contrast to this, if $ \{u\}_0 $ is a set of ordered $ \forall u_i \in [c,d] $ the superset $ \{u\}_0 \supseteq O $ has only $ W=1 $ elements and its entropy is $ S_{\{O\}}=0 $.  Shannon \citep{Szilard1929, Shannon1948} information is $  \iota = -k \, \ln(W) $, with $ k $ a constant which may may be taken as $ k=1 $, is the negative of entropy (or \textit{negentropy}); if we know that an unknown message contains $ m $ words we face $ m! $ possibilities, the information we have is $ \iota= k \ln (m!) $, yet when the message is decoded we are left with a single possibility and $ \iota = k \ln(1)=0 $ the maximum negentropy, or the maximum information possible by Shannon\textquoteright{s} \citep{Shannon1948} definition. An important particular case occurs when $c=0 $, $ d=1 $ and the uniform distribution is defined for the closed unit interval $ [0,1] \in \mathbb{R}$ denoted here as $ U[0,1], $ which enables Monte Carlo statistical simulation, described in more detail ahead. Central to this paper is the rectangular distribution where only the digits of a number in base $ b $ are considered, and are equally likely with probability $ b^{-1} $. An infinite sequence of digits fulfilling Eq. (\ref{eq:IntSec}) is distributed as  $ U_{\mathbb{Z}}[0, b-1] $, has a negentropy $ \iota=- \infty $, and entropy $ S_\Omega=\infty $.

In this paper the digits of $\pi$ are considered as a sequence (called here also chain, catena, catenation succession or string to reduce repetitiveness) sampled from a waveform at constant intervals, and the properties of such waveform are studied using fractal analysis. The fractal dimension may convey information on spatial extent (convolutedness or space filling properties), self-similarity (the ability to remain unchanged when the scale of measurement is changed) and self affinity \citep{Barnsley1993}. In signal (waveform) digital processing the continuous waveform is sampled at regular intervals; this sampling results in set $ \{y_i \}_{i=1,2,\dots,N}$ which is further processed numerically to extract its information content \citep{Smith1997}.  

An expression to calculate the fractal dimension of a waveform was obtained by  \citet{Sevcik1998a} starting from the definition of Hausdorff-Besicovitch \citep{Hausdorff1918, Besicovitch1929} dimension ($D_h$). Mandelbrot\textquoteright{s} definition calls fractal \citep{Mandelbrot1983} to a set whose Hausdorff-Besicovitch dimension is not an integer. The Hausdorff-Besicovitch dimension of a set in a metric space  may be expressed as:
\begin{equation}\label{eq:DH}
D_{h}=-\underset{\varepsilon \rightarrow 0}{\lim }{\frac{\ln [N(\varepsilon )]}{\ln (\varepsilon)}}
\end{equation}
where $N(\varepsilon)$ is the number of open balls of a radius $ \varepsilon $ needed to cover the set. Given a point $ p $ and defining distance between $p$ and another point $x$ in the same space as $ \delta(x,p) $, an open ball of center $ p $ and radius $ \varepsilon $, is a set of all points $ x $ for which $ \delta(p,x) < \varepsilon$.

Waveforms are planar curves in a space with coordinates usually having different units. Since the topology of a metric space does not change under linear transformation, it is convenient linearly to transform a waveform into another in a normalized space, where all axes are equal. \citet{Sevcik1998a} proposed the use of two linear transformations mapping the original waveform into another embedded in an equivalent metric space. The first transformation, normalizes every point in the abscissa as:
\begin{equation}\label{eq:xtrans}
x_i^*  = \frac{x_i-x_{min}}{x_{max}-x_{min}}
\end{equation}
where $x_i $ are the original values of the abscissa, and $x_{max}$ and $x_{min}$ are the minimum and maximum $x_i$, respectively. The second transformation normalizes the ordinate as follows:
\begin{equation}\label{eq:ytrans}
y_i^*  = \frac{y_i-y_{min}}{y_{max}-y_{min}}
\end{equation}
where $y_i $ are the original values of the ordinate, and $y_{max}$ and $y_{min}$ are the minimum and maximum $y_i$, respectively. These two linear transformations map the $N$ points of the waveform into another that belongs to a unit square. This unit square may be visualized as covered by a grid of $N \cdot N$ cells. $N$ of them containing one point of the transformed waveform. Calculating the length $L$ of the embedded waveform, the length of the transformed waveform, and taking $\epsilon=\frac{1}{2 \cdot N'}$ Eq. (\ref{eq:DH}) becomes \citep[Eq. (6a)]{Sevcik1998a}
\begin{equation}\label{eq:DCS}
D_h = \Phi \approx D = 1 + \underset{N' \rightarrow \infty}\lim \left[ \dfrac{\ln(L/2)}{\ln(2 \cdot N')}\right] .
\end{equation}

If $ \pi $ is normal, then its digits must be a  sequence fulfilling the condition expressed by Eq. (\ref{eq:IntSec}) and distribured as $ U_{\mathbb{Z}}[0,9] $. Considering the $ \pi $ digit catenation ¿as a set $\{y_i\}_{i=0,1,\ldots,N-1}  $ of points sampled at constant $ \Delta x = x_{i+1}-x_i $, from an hypothetical waveform approximated by a set of straight segments, each one of them with length equal to
 \[ \Delta L_{i}= \sqrt{(y_{i+1}-y_i)^2 +\Delta x^2 }\]
 which are transformable as indicated by Eqs. \ref{eq:xtrans} and  \ref{eq:ytrans} into the embedded form
 \begin{equation}
 \Delta L^*_{i}= \sqrt{(y^*_{i+1}-y^*_i)^2 +(\Delta {x^*})^2 }
 \end{equation}
 which adds
 \begin{equation}
 L^*= \sum_{i=0}^{N'}\Delta L^*_{i}\text{,}
 \end{equation}
 where $ L^* $ is the digital approximation of total transformed waveform length. Eq. (\ref{eq:DCS}) may now be used to estmate $ D $ for the $ N' $ points. As indicated in Section \ref{sec:StatCons} and the condition expressed by Eq. (\ref{eq:Lim2}) $ \underset{N' \rightarrow{} \infty}{\lim} D= \Phi $, but will converge through a different pathway for different waveforms with rhe same $ \Phi $ or to different values of $ \Phi $ in other waveforms \citep{Sevcik1998a}. When the limiting value and the manner it is reached are considered together, they become a fingerprint to differentiate between waveforms.

\section{Methods.}\label{sec:MatMet}

\subsection{Generating 10\textsuperscript{9} digits of {\greektext{p}}.}

One billion decimal digits of $\pi$ were generated [Total time (base 10 result) = 1929 s] using the tpi-0.9 executable multiprocessor threaded program for Linux provided by \citet{Bellard2010}, and were stored on on disk for further analysis. The tpi-0.9 code implements an improved Ramanujan series \citep{Chudnovsky2000} and has been used to successfully produce up to $2.7\!\cdot\!10^{12} \; \pi$ digits under similar conditions \citep{Bellard2010}. The limit of $ 10^9 $ digits was decided due to restrictions in the computer power available to the author. 

\subsection{Random Number Generation.}\label{sec:Random}

Fundamental to all Monte Carlo simulations \citep{Dahlquist1974} is a good uniform (pseudo) random (PRNG) number generator. Data for all numerical simulations carried out in this work were produced using  random numbers ($r$) with continuous rectangular (uniform) distribution in the closed interval [0,1]. To avoid singularities in the simulations of exponentially distributed pseudo random numbers, $r \in (0,1], \;(U(0,1]$), were used. All $U[0,1]$  or $U(0,1]$ were generated using the 2002/2/10 initialization-improved  623-dimensionally equidistributed uniform pseudo random number generator, MT19937 algorithm \citep{Matsumoto1998, Panneton2006}. The generator has passed the stringent DIEHARD statistical tests \citep{Marsaglia2003, Bellamy2013}. It uses 624 words of state per generator and is comparable in  speed to the other generators.  It has a Mersenne prime period of $2^{19937} -1$ ($\approx10^{6000}$).  

The MT19937 seed used was a 64-bit unsigned integer obtained using the /dev/random Linux PRNG, which employs environmental noise from device drivers and other sources into an entropy pool. Device /dev/random gets blocked, and stops producing random bytes, if the entropy of the device gets low, and commences producing output again when it recovers to safe levels. No such delays were perceived during this work. Using /dev/random seed makes exceedingly unlikely ($P = 2^{-64} \approx 5.4 \cdot 10^{-20}$) that the same sequence of $U[0,1]$, $\{r_i\}$, is used twice.

The following types of $\{r_i\}$ were generated and used in this work:
\begin{enumerate}
\item Uniform random decimal digits $ U_{\mathbb{Z}}[0,9] $ were generated as the integer part of $10 \cdot U[0,1)$. The integer part was obtained with the int() typecasting in C++ \citep{Stroustrup2013}.
\item Random normal Gaussian variates with mean 0 and variance 1, $N(\mu =0,\sigma^2=1)\in\mathbb{R}$, denoted from now on as $N(0,1)$, were generated with the \citet{Box1958} algorithm implemented in C++ as a variation of the gasdev() C function  in \citet[pg. 293]{Press1992}.
\item Random Poisson variables  \[ Po(x | \eta)=\frac{\eta^x}{x !}e^{-\eta} \]  with $\mu=\sigma^2=1$, $Po(x | \eta=1)$, were generated for $x \in  \mathbb{Z}[0,9]$.
\item Brownian random walks of length $N$ were generated recursively as 
\begin{equation}\label{eq:BrownW}
\lbrace r_{i+1} \looparrowleft r_i + N(0,1) \rbrace_{i=2, \ldots ,N-1}
\end{equation}

starting at $r_1=N(0,1)$.
\item Uniform random walks of length $N$ were generated recursively as 
\begin{equation}\label{eq:UnifW}
\lbrace r_{i+1} \looparrowleft r_i +(2 \cdot U[0,1]-1) \rbrace_{i=2, \ldots ,N-1}
\end{equation}
starting at $r_1=2 \cdot U[0,1] -1$.
\end{enumerate}

To test randomness of some digit strings, $ D $ values of these sequences were compared before and after randomization (see Tables \ref{tab:PiPiRand} and  \ref{tab:RandPervsPer}). In these instances, randomization was achieved by generating a vector of $N=10^9$ real numbers ($\{r_i\}_{i=1,2,\ldots,N}$) in $U[0,1]$ which were sorted with the heapsort algorithm \citep{Williams1964, Schaffer1993}. The heapsort algorithm was coded to transpose the $d_i$ digit  of $\pi$ ($\{d_i\}_{i=1,2,\ldots,N}$) at the same time as the isomorphic  $r_i$ was transposed in the sorting process, producing a totally randomized set $\{d_i\}$.. 

All computations in this paper, were programed in C++14 \citep{Stroustrup2013} and compiled with the GNU gcc/g++ compiler (version 5.2.1 20151010, Ubuntu 5.2.1-22ubuntu2, htpp://gcc.gnu.org) under 64 bit Ubuntu Linux version 15.10 on an Apple MacBook Air computer (8 GB RAM, Intel\textsuperscript{\textregistered} Core\texttrademark{} i7-4650U CPU @ 1.70 GHz $ \times $ 4, RAM disk 500 GB).  Source codes of all programs used in this work are included on line as supplementary material.

\subsection{Statistical Considerations}\label{sec:StatCons}

\subsubsection{$D$ is convergent as $N'$ increases.}\label{sec:Converg}

As indicated in \citet{Sevcik1998a}, although the fractal dimension $\Phi$ is a topological invariant of a set or a metric space, \textit{D} is only an empirical estimate of $\Phi$ with some uncertainty based on a set of points sampled from a waveform; \textit{D} is thus a random variable. The relationship between $\Phi$ and \textit{D}, is similar to the one between the mean of a population ($\mu$) and the mean $\overline{x}$ estimated sampling a subset of the population; although $ \mu $ is an invariant for the population, $\overline{x}$ will change with sampling. Just as $\overline{x}$ converges to $ \mu $ as the sample size approaches the size of the population, \textit{D} converges to $\Phi$ as $N' \longrightarrow \infty$.The expression for the variance of $D$ was determined \citep[Eqs. (10) and (11)]{Sevcik1998a} as 

\begin{equation} \label{eq:VarDemp}
\Var(D)=\frac{N'\cdot \Var(\Delta y)}{L^{2}\cdot \ln (2\cdot N')^{2}} 
\end{equation}
where $\Var (\Delta y)$ may be estimated from the data as:
\begin{equation}\label{eq:VarDeltay}
\Var(\Delta y)=\overset{N'}{\underset{i=1}{\sum }}{\frac{\left(\Delta y_{i}-{\overline{{\Delta y}}}\right)^{2}}{N'}}
\end{equation}
where  ${\overline{{\Delta y}}}$ is the mean segment length. The limit of sample mean variance is
\begin{equation}
\underset{N' \rightarrow \infty}\lim \Var( \Delta y) = \sigma^2( \Delta y)=2\, \sigma^2(y)
\end{equation}
in which $\sigma^2( y)$ is the variance of the function used to generate the sequences of variables described in Section \ref{sec:Random}. Thus.
\begin{equation}\label{eq:VarD}
\Var(D)=\frac{2 \cdot N' \cdot \sigma^2(y)}{L^{2}\cdot \ln (2\cdot N')^{2}}
\end{equation}
and
\begin{equation}\label{eq:LimVarD}
\underset{N' \rightarrow \infty}\lim \Var( D) =\sigma^2(D) = \sigma^2(y) \cdot \underset{N' \rightarrow \infty}\lim \left \{  \frac{2 \cdot N' }{L^{2}\cdot \ln (2\cdot N')^{2}} \right \}.
\end{equation}
The length of the embedded waveform is a sum of straight line segments joining the sample points in the waveform: 
\begin{equation}
L=\overset{N'}{\underset{i=1}{\sum }} \sqrt{\Delta y_i^2+N'^{-2}}
\end{equation}
for very large N'
\begin{equation}
\begin{split}
\underset{ N'\rightarrow \infty}\lim L^2 = \underset{N' \rightarrow \infty}\lim\left\lbrace  \overset{N'}{\underset{i=1}{\sum }} \sqrt{(N' \overline{\Delta y_i})^2+N'^{-2}}\right\rbrace^2\\
\cdots = \underset{N' \rightarrow \infty}\lim\left\lbrace \overset{N'}{\underset{i=1}{\sum }} N' \overline{\Delta y_i}\right\rbrace^2 &= \underset{N' \rightarrow \infty}\lim  N'^2 \left\lbrace {\overset{N'}{\underset{i=1}{\sum }}  \overline{\Delta y_i}}\right\rbrace ^2
\end{split}
\end{equation}
thus
 \begin{equation}\label{eq:ConvVarD}
 \underset{N' \rightarrow \infty} \lim [\Var(D)\;|\;0<\sigma^2(y)<\infty]=0. 
 \end{equation}
 Thus for very large $N'$, $D$ converges to a constant value.
 
\subsection{Statistical considerations on Monte Carlo \textit{D} data analyses}
\subsubsection{The Probability Density Function of $D$.}

Crucial to this study is to determine whether $ D $ values calculated under different conditions are due to the difference in condition or whether the difference stems from sampling variation and uncertainty. In general the probability density function of \textit{D} is unknown, but there is empirical evidence that is  unimodal but not Gaussian and is strongly skewed to the right (Skewness $  > 0$  $ $); see for example Figure 3 in \citet{DSuze2010a}. It is tempting to assume that with $ N' $ points large, the distribution of  ${\Delta{y}}$ in Eq.~(\ref{eq:VarDeltay}) is asymptotically Gaussian as demonstrated by the central limit theorem \citep{Wilks1962}. However, infinity is far away for some pdfs and thus Gaussianity must be demonstrated to avoid statistical errors if parametric test are to be used. All test for Gaussianity have shortcomings. Two of the most powerful and modern test to determine whether a set of random variables come from a population having a Gauss pdf are the \citet{Bera1981} and  \citet{Shapiro1965} tests, thus both test were used to assert Gaussianity of $\{D_i\}$ sets compared in this paper. Results of the tests were quite variable, suggesting that some sets were highly non-Gaussian ($P<10^{-6}$) while in most small samples $0.1<P<0.5$, perhaps reflecting low power of the tests with sample sizes below 100. Due to this we present both nonparametric and parametric analysis results. We prefer the more conservative nonparametric test and the \citet{Vysochanskij1980} inequality (Theorem \ref{th:VP}) to decide against the null hypothesis since this ensures that differences declared significant with the nonparametric test are certainly more statistically significant than they appear.

\subsubsection{On nonparametric test used.}\label{sec:nonparam}
 
In this paper $ \{D\} $ sets with uncertainty are compared under different Monte Carlo statistical experiment conditions. Parametric comparison with Gaussian tests are generally more powerful than their non-Gaussian alternatives if and only if the data pdf is $N(\mu,\sigma^2)$ but may result in statistical errors of type I or II if the pdf is not Gaussian. When the pdf is not Gaussian, nonparametric alternatives are more powerful, less prone  to produce errors. Thus, all conclusions from the Monte Carlo simulations in this paper were checked for statistical significance with nonparametric tests and the Vysochanskij-Petunin inequality (Section \ref{sec:VP}). 
  Multiple comparisons were done with the  \citet{Kruskal1952} nonparametric analysis of the variance. Significances based on parametric test are included just for comparison purposes, not for decision making. For all nonparametric test used but not referenced to primary sources please check \citet{Hollander1973}. Although significances are usually expressed in the conventional manner $ P=\alpha $, meaning that there is $ \alpha $ probability that a difference between tests of random variables stems from chance, on some occasions the notation $ P(\text{H0}) = \alpha $ is alternatively used to indicate the probability of the null hypothesis (H0).

\subsubsection{The Vysochanskij-Petunin Inequality.}\label{sec:VP}

The probability that a random variable $x$ comes from a Gaussian population with mean $\mu$ and variance $\sigma^2$ [$N(\mu,\sigma^2)$] may be estimated defining
\begin{equation}\label{eq:z}
\lambda=\frac{ \abs{x-\mu} }{\sigma}
\end{equation}
and calculating
\begin{equation}\label{eq:Gauss}
\varPsi(\lambda) = 1- \varPhi(\lambda)=1-\frac{1}{\sqrt{2 \pi}}\int\limits_{-\infty}^{\lambda} \mathrm{e}^{-\frac{1}{2}z^2}\mathrm{d}z,
\end{equation}
where $\varPhi(\lambda)$ is the Gauss probability distribution function (\textit{PDF}) \cite{Wilks1962}. When $\varPsi(\lambda)\leq \epsilon$ the null hypothesis, $x \in N(\mu,\sigma^2)$, is rejected with $\geq(1-\epsilon)$ certainty, the so called $P\leq\epsilon$ \textit{confidence level}. 

If the pdf of the random variable is not known, an answer may still be obtained with the Vysochanskij-Petunin inequality \citep{Vysochanskij1980}. The Vysochanskij-Petunin inequality provides an upper bound for the probability that a random variable with finite variance lies within a certain number of standard deviations from the variable's mean; or equivalently, it provides an upper bound for the probability that it lies further away. The sole restrictions are that the distribution is unimodal and has finite variance. The inequality requires a continuous probability distribution, except perhaps at the mode which may have a non-zero probability. 

\begin{theorem}[Vysochanskij-Petunin \citep{Vysochanskij1980}]\label{th:VP}
Let $x$ be a random variable with unimodal distribution, mean ${\mu}$ and finite, non-zero variance ${\sigma^2}$. For any \[{\lambda > \sqrt{\frac{8}{3}} \approx 1.63299\ldots}\] then
\begin{equation}\label{eq:VP}
    P(\left| x- \mu \right| \geq \lambda \sigma) \leq \frac{4}{9\lambda^2}=\epsilon.
\end{equation}
\end{theorem}

The theorem applies even to heavily skewed distributions and puts bounds on how much of the data is, or is not, \textquotedblleft{}in the middle\textquotedblright. Setting ${\alpha = 0.05}$ then ${\lambda=\pm \left ( \sqrt{80/9}\approx 2.981424 \ldots \right )}$, for a two tailed test. By virtue of Eq.~(\ref{eq:VP}) no matter which unimodal distribution, no matter how skewed, there will be $<2.5$\% chance that a datum will belong to a population with mean $\mu$ and variance $\sigma^2$ if it lays farther than $\pm2.981\ldots\:\sigma$ from $\mu$. Please note that if the probability distribution function of data is Gaussian $\pm1.96\ldots\sigma$ suffices to reach the same confidence level. The $ \epsilon $ variable in Eq.~(\ref{eq:VP}) was introduced by the author of the present communication for the statistical argumentation which follows.

Let $\overline{x_1}$ and  $s^2(x_1)$ be the mean and variance estimated for the random variable ${x_1}$, and $\overline{x_2}$ and  $s^2(x_2)$ be the mean and variance estimated for the random variable ${x_2} $ independent of $ {x_1}$, then $\Delta{\overline{x}_{1,2}}=\overline{x_1}-\overline{x_2}$. The variance of $\Delta{\overline{x}_{1,2}}$ is $\Var(\Delta{\overline{x}_{1,2}})=\Var(x_1)+\Var(x_2)$. By virtue of the Vysochanskij-Petunin inequality~[Eq. (\ref{eq:VP})]

\begin{equation}\label{eq:PDif}
P \left(\frac{ \left|\Delta{\overline{x}_{1,2}} \right|}{s(\Delta{\overline{x}_{1,2}})} \geq \sqrt{\frac{4}{9 \epsilon}}\right ) \leq \epsilon
\end{equation}
and thus
\begin{equation}\label{eq:PVP}
\epsilon = \frac{2}{9 \left[ \frac{\Delta{{\overline{x}}_{1,2}}}{s(\Delta{\overline{x}_{1,2}})}\right ]^2}\text{ .}
\end{equation}

Therefore, there is $\leq\epsilon$ probability that $\left|\Delta{\overline{x}_{1,2}} \right|\neq0$ due to random sampling variation, and $ \Delta{\overline{x} }_{1,2} \neq 0 $ with a confidence level $P \leq \alpha = \epsilon$. Due to this
\begin{equation}\label{eq:P_VP2}
P \left[\frac{ \left|\Delta{\overline{x}_{1,2}} \right|}{s(\Delta{\overline{x}}_{1,2})} \geq \left ( \sqrt{\frac{40}{9}}\approx2.108\ldots\right )\right ] \leq 0.05
\end{equation}
for a test with only one tail, because only one alternative matters: $\left|\Delta{\overline{x}_{1,2}} \right|\neq0$. When the probability distribution function is Gaussian, the value of $\lambda$ for the one tailed case, would be $\approx1.65\ldots$ instead of $\approx2.108\ldots$ demanded by Eq.~(\ref{eq:VP}). 
\begin{lemma}\label{th:alpha}
When the conditions of Theorem \ref{th:VP} are fulfilled, and $ \alpha $ is our decision level, then $P(\Delta{\overline{x}_{1,2}}\neq 0\, |\, \epsilon \geqq \alpha) \leqq \alpha $ and $ P(\Delta{\overline{x}_{1,2}}\neq0\, |\, \epsilon \leq \alpha) < \alpha $, due to random sampling variation.
\end{lemma}

\section{Results.}\label{sec:Results}

\begin{table}[h!]
\begin{center}
\caption{{Evaluating $ D $ for the $ \pi $ digits as $ N' $ grows from 10 to $ 10^{9} $ before and after randomizing the sequence}}
\begin{tabular}{*{1}{m{3cm}}*{1}{m{3.3cm}}*{1}{m{1cm}}*{1}{m{1.3cm}}c}
\hhline{*{5}{=}}
\centering $N'$ &\centering $ \overline{D} \pm sd $ & \centering $ \lambda$ & \centering $ \epsilon $& $ \varPsi(\lambda) $\\
\hhline{*{5}{-}}
\centering$ \pi $&&&&\\
\centering 10&\centering 1.229738 \textpm{} 0.055902&&&\\
\centering 100&\centering 1.558271 \textpm{} 0.012013&&&\\
\centering 1,000&\centering 1.689786 \textpm{} 0.002944&&&\\
\centering 10,000&\centering 1.758991 \textpm{} 0.000721&&&\\
\centering 100,000&\centering 1.804246 \textpm{} 0.000186&&&\\
\centering 1,000,000&\centering 1.835332 \textpm{} 0.000049&&&\\
\centering 10,000,000&\centering 1.857853 \textpm{} 0.000014&&&\\
\centering 100,000,000&\centering 1.874979 \textpm{} 0.000004&&&\\
\centering 1,000,000,000&\centering 1.888421 \textpm{} 0.000001&&&\\
\centering{Randomized $ \pi $}&&&&\\
\centering 10&\centering 1.239781 \textpm{} 0.071590&0.1106&\centering NA&0.46\\
\centering 199&\centering 1.561186 \textpm{} 0.013844&0.1590&\centering NA&0.44\\
\centering 1,000&\centering 1.684466 \textpm{} 0.003061&1.2530&\centering NA&0.11\\
\centering 10,000&\centering 1.757939 \textpm{} 0.000738&1.0195&\centering NA&0.15\\
\centering 100,000&\centering 1.804153 \textpm{} 0.000186&0.3517 &\centering NA&0.3 \\
\centering 1,000,000&\centering 1.835291 \textpm{} 0.000049&0.5891 &\centering NA&0.28\\
\centering 10,000,000&\centering 1.857851 \textpm{} 0.000014&0.1285  &\centering NA &0.45\\
\centering 100,000,000&\centering 1.874980 \textpm{} 0.000004 &0.2138 &\centering NA &0.41\\
\centering 1,000,000,000&\centering 1.888422 \textpm{} 0.000001 &0.6154 &\centering NA&0.27\\
\hhline{*{5}{=}}
\end{tabular}\label{tab:PiPiRand}
\end{center}
\footnotesize{ $D$: Calculated with Eq. \ref{eq:DCS} ; sd: standard deviation of $D$ calculated with Eq. \ref{eq:VarDemp} ; $N' $: number of $ \pi $ digits used to calculate $D$; {\selectlanguage{greek}l}: Calculated as in Eqs. \ref{eq:VP} and \ref{eq:PDif} for the difference between $D$ before and after of randomizing $ \pi $; $ \varPsi(\lambda) $: Probability calculated with the Gauss pdf (Eq. \ref{eq:Gauss}) for the difference in $D$ before and after of randomizing {\selectlanguage{greek}p}; NA: Not applicable.}
\end{table}

\subsection{Fractal analysis of {\greektext{p}} digits sequence.} \label{sec:FracAnPi}

Table \ref{tab:PiPiRand} presents $D$ for the $\pi$ digits series before (Top of the table) and after randomizing the sequence as indicated in Section \ref{sec:Random} (Bottom of the table). In both halves of the table, values of $D$ were obtained initially with the first 10 digits ($ N'=10 $), and subsequently increasing $ N' $ in steps of ten times the previous \textit{N'} value. Data in Table \ref{tab:PiPiRand} shows, that $D$ becomes less variable and approaches 2 as $N'$ increases, and that this is true for both native and randomized $ \pi $ digit sequences. 

The values  $ D $ of the $ \pi $ digits sequence are compared in Table \ref{tab:PiPiRand} prior to and after randomizing this sequence; for this ourpose values of $\lambda$, $ \varepsilon $ (Eq. \ref{eq:PVP}) and $ \varPsi(\lambda) $ (Eq. \ref{eq:Gauss}) are presented at right of the bottom half of the table. $ \varPsi(\lambda) $ is the probability calculated with the Gauss pdf (Eq. \ref{eq:Gauss}) for the difference in $D$, calculated with the same $ N' $, before and after randomizing the $ \pi $ sequence. All values of $ \lambda$ in the table were  $< \sqrt{8/3} $ which means that $\Delta D$ is too small to be statistically significant if the Vysochanskij-Petunin inequity is used to compare (NA in the table), and furthermore that they are too small for the  Vysochanskij-Petunin theorem to hold. The values of $ \varPsi(\lambda) $, calculated assuming Gaussianity of data, also indicates that $\Delta D \neq 0$, as it should be if randomization does not chance $ \pi $ sequence\textquoteright{s} $ D $. The data also show that the more $\pi$ digits are considered, the more the sequence resembles a concatenation of uniform random independent variables $U_{ \mathbb{Z}} [0,9] $.

An analysis and notation entirely similar are used in the following tables.
\begin{table}[h!]
\begin{center}
\caption{{Comparing $ D $ for the $ \pi $ digits as $ N' $ grows from 10 to $ 10^{9} $ with $ \overline{D} $  calculated for 30 sequences of random independent and uniformly distributed ($U_{\mathbb{Z}}[0,9]$) decimal integers and $ U[0,1] $ real numbers under similar conditions.}}
\begin{tabular}{*{1}{m{3cm}}*{1}{m{3.3cm}}*{1}{m{1cm}}*{1}{m{1.3cm}}c}
\hhline{*{5}{=}}
\centering $N'$ &\centering $ \overline{D} \pm sd $ & \centering $ \lambda$ & \centering $ \epsilon $& $ \varPsi(\lambda) $\\
\hhline{*{5}{-}}

\centering $ U_{\mathbb{Z}}[0,9]$ Uniform &&&\\
\centering 10&\centering 1.235588 \textpm{} 0.069559&0.084& \centering NA&0.46\\
\centering 100&\centering 1.550211 \textpm{} 0.018428&0.437& \centering NA&0.33\\
\centering 1,000&\centering 1.684800 \textpm{} 0.003208&1.554& \centering NA&0.06\\
\centering 10,000&\centering 1.758777 \textpm{} 0.000871&0.246& \centering NA&0.40\\
\centering 100,000&\centering 1.804223 \textpm{} 0.000213&0.109& \centering NA&0.46\\
\centering 1,000,000&\centering 1.835294 \textpm{} 0.000046&0.835& \centering NA&0.20\\
\centering 10,000,000&\centering 1.857650 \textpm{} 0.001141&0.178& \centering NA&0.43\\
\centering 100,000,000&\centering 1.874429 \textpm{} 0.001682&0.327& \centering NA&0.37\\
\centering 1,000,000,000&\centering 1.887438 \textpm{} 0.002002&0.491& \centering NA&0.31\\
\centering$ U[0,1]$ Uniform &&&\\
\centering 10&\centering 1.190863 \textpm{} 0.073295 &0.530& \centering NA&0.298\\
\centering 100&\centering 1.530226 \textpm{} 0.012791&2.193& \centering 0.092 &0.014\\
\centering 1,000 &\centering 1.671170 \textpm{} 0.003572&5.211& \centering 0.016 &$9 \cdot 10^{-8}$\\
\centering 10,000&\centering 1.748951 \textpm{} 0.000740 &13.568& \centering 0.002 &$<10^{-14}$\\
\centering 100,000&\centering 1.796415 \textpm{} 0.000227&34.46&\centering $4\cdot 10^{-4}$ &$<10^{-14}$\\
\centering 1,000,000&\centering 1.828722 \textpm{} 0.000065 &101.5&\centering $4\cdot 10^{-5}$ &$<10^{-14}$\\
\centering 10,000,000&\centering 1.852191 \textpm{} 0.000012 &459.7&\centering $2\cdot 10^{-6}$ &$\ll10^{-14}$\\
\centering 100,000,000&\centering 1.869994 \textpm{} 0.000003 &1630&\centering $2\cdot 10^{-7}$ &$\ll10^{-14}$\\
\centering 1,000,000,000&\centering 1.883972 \textpm{} 0.000001 &4828&\centering $2\cdot 10^{-8}$ &$\ll10^{-14}$\\

\hhline{*{5}{=}}
\end{tabular}\label{tab:PiVsU09}
\end{center}
\footnotesize{ $D$: Calculated with Eq. \ref{eq:DCS} ; sd: standard deviation of $D$ calculated with Eq. \ref{eq:VarDemp} ; $N' $: number of $ \pi $ digits used to calculate $D$; $ \mathbb{Z}[0,9] $calculated as indicated in Section \ref{sec:Random} ; $ \lambda $: Calculated as in Eqs. \ref{eq:VP} and \ref{eq:PDif} for the difference between $D$ before and after of randomizing $ \pi $; $ \varPsi(\lambda) $: Probability calculated with the Gauss pdf (Eq. \ref{eq:Gauss}) for the difference in $D$ before and after of randomizing {\selectlanguage{greek}p}; NA: Not applicable. Values for $\pi$ in Table \ref{tab:PiPiRand}.}
\end{table}
\subsection{Comparing the fractal properties of {\greektext{p}} sequence with uniformly distributed sequences.} 

The $ \pi $ sequence was compared (Top of Table \ref{tab:PiVsU09}) with the average $D$ obtained for 30 sequences of random uniform and independent decimal digits ($U_{\mathbb{Z}} [0,9] $) generated as indicated in Section \ref{sec:Random}. The standard deviations presented for series in the top of Table \ref{tab:PiVsU09} include the waveform sampling variance expressed by Eq. \ref{eq:VarDemp} but also the variability between the 30 sets of $U_{\mathbb{Z}} [0,9]  $ generated. As seen in the top of Table \ref{tab:PiVsU09}, all $ \lambda $ values are too small to detect a significant difference between the $ D $s calculated for $ \pi $ and for $  U_{\mathbb{Z}}[0,9] $ sets. This suggests that $ \pi $ digits sequence is of type $ U_{\mathbb{Z}}[0,9] $.

The bottom of Table \ref{tab:PiVsU09} is similar to its top, except for that a sequence of real numbers uniformly distributed in the closed interval $ [0,1] $ ($ U[0,1] $) is compared with the $ \pi $ sequence. As seen in the table, all $ U[0,1] $ real number sequences with $ N'>100 $ were statistically different in  $ D $ from $ \pi $ digit series ($ P \leqq 0.05 $), which shows that fractal analysis [Eq. (\ref{eq:DCS})] is useful to distinguish between random uniformly distributed and $\pi $ digit  chain, and confirms that $ \pi $ digits do not form a $ U[0,1] $ sequence.

\begin{table}[h!]
\begin{center}
\caption{{Comparing $ D $ for the $ \pi $ digits as $ N' $ grows from 10 to $ 10^{9} $ with $ \overline{D} $ calculated for 30 sequences of random independent Gaussian ($N(0,1$), exponential [$ f_e(x|\eta=1) $] or Poissonian [$f_p(x | \eta=1)$] variables.}}
\begin{tabular}{*{1}{m{3cm}}*{1}{m{3.3cm}}*{1}{m{1cm}}*{1}{m{1.3cm}}c}
\hhline{*{5}{=}}
\centering $ N'$ &\centering $ \overline{D} \pm sd $ & \centering $ \lambda$ & \centering $ \epsilon $& $ \varPsi(\lambda) $\\
\hhline{*{5}{-}}

\centering Gaussians &&&\\
\centering 10&\centering 1.192427 \textpm{} 0.070534&0.529& \centering NA&0.298\\
\centering 100&\centering 1.465974 \textpm{} 0.024732&3.731& \centering 0.032&$3\cdot10^{-5}$\\
\centering 1,000&\centering 1.590481 \textpm{} 0.010237&9.701& \centering $5\cdot 10^{-3}$ &$<10^{-14}$\\
\centering 10,000&\centering 1.665459 \textpm{} 0.006783&13.79& \centering $2\cdot 10^{-3}$ &$<10^{-14}$\\
\centering 100,000&\centering 1.719042 \textpm{} 0.003256&26.17&\centering $6\cdot 10^{-4}$ &$<10^{-14}$\\
\centering 1,000,000&\centering 1.756622 \textpm{} 0.003274&24.04&\centering $8\cdot 10^{-4}$ &$<10^{-14}$\\
\centering 10,000,000&\centering 1.784259 \textpm{} 0.001713&42.97&\centering $2\cdot 10^{-4}$ &$\ll10^{-14}$\\
\centering 100,000,000&\centering 1.806736 \textpm{} 0.001254&54.40&\centering $1.5\cdot 10^{-4}$ &$\ll10^{-14}$\\
\centering 1,000,000,000&\centering 1.824420 \textpm{} 0.000874&73.23&\centering $8\cdot 10^{-5}$ &$\ll10^{-6}$\\
\centering Exponential &&&\\
\centering 10&\centering 1.174037 \textpm{} 0.085295&0.653& \centering NA&0.257\\
\centering 100&\centering 1.435551 \textpm{} 0.032232&3.807& \centering 0.031 &$7\cdot 10^{-5}$\\
\centering 1,000&\centering 1.555341 \textpm{} 0.017704&7.594& \centering $8\cdot 10^{-4}$ &$2 \cdot 10^{-14}$\\
\centering 10,000&\centering 1.629120 \textpm{} 0.011993&10.83& \centering $3\cdot 10^{-3}$ &$<10^{-14}$\\
\centering 100,000&\centering 1.683890 \textpm{} 0.007756&15.52&\centering $2\cdot10^{-3}$ &$<10^{-14}$\\
\centering 1,000,000&\centering 1.721709 \textpm{}0.006243&18.20&\centering $1\cdot 10^{-3}$ &$<10^{-14}$\\
\centering 10,000,000&\centering 1.750673 \textpm{} 0.003643&29.42&\centering $5\cdot 10^{-4}$ &$\ll10^{-14}$\\
\centering 100,000,000&\centering 1.774054 \textpm{}0.003211&54.40&\centering $4\cdot 10^{-4}$ &$\ll10^{-14}$\\
\centering 1,000,000,000&\centering 1.792871 \textpm{} 0.002229&42.87&\centering $2\cdot 10^{-4}$ &$\ll10^{-14}$\\
\centering Poissonian &&&\\
\centering 10&\centering 1.252291 \textpm{} 0.107372&0.210& \centering NA&0.417\\
\centering 100&\centering 1.555190 \textpm{} 0.022202&0.139& \centering NA &0.445\\
\centering 1,000&\centering 1.690384 \textpm{} 0.004576&0.131& \centering NA &0.448\\
\centering 10,000&\centering 1.762284 \textpm{} 0.001184&2.783& \centering  0.057&$3\cdot 10^{-3}$\\
\centering 100,000&\centering 1.807332 \textpm{} 0.000339&9.12&\centering $5\cdot 10^{-3}$ &$<10^{-14}$\\
\centering 1,000,000&\centering 1.837881 \textpm{} 0.000093&27.45&\centering $6\cdot 10^{-3}$ &$\ll10^{-14}$\\
\centering 10,000,000&\centering 1.860074 \textpm{} 0.000025&90.36&\centering $5\cdot 10^{-5}$ &$\ll10^{-14}$\\
\centering 100,000,000&\centering 1.876927 \textpm{} 0.000006&315.5&\centering $5\cdot 10^{-6}$ &$\ll10^{-14}$\\
\centering 1,000,000,000&\centering 1.890159 \textpm{} 0.000002&1058&\centering $4\cdot 10^{-7}$ &$\ll10^{-14}$\\
\hhline{*{5}{=}}
\end{tabular}\label{tab:GausPoiExp}
\end{center}
\footnotesize{{$D$: Calculated with Eq. \ref{eq:DCS} ; sd: standard deviation of $D$ calculated with Eq. \ref{eq:VarDemp} ; $N' $: number of $ \pi $ digits used to calculate $D$; $ U_{\mathbb{Z}}[0,9] $calculated as indicated in Section \ref{sec:Random} ; $ \lambda $: Calculated as in Eqs. \ref{eq:VP} and \ref{eq:PDif} for the difference between $D$ before and after of randomizing $ \pi $; $ \epsilon $: Probability calculated with the Vysochanskij-Petunin inequality (Eq. \ref{eq:PVP}) for the difference in $D$ before and after of randomizing {\selectlanguage{greek}p}; NA: Not applicable. Values for $\pi$ in Table \ref{tab:PiPiRand}.}}
\end{table}

\subsection{Comparing the fractal properties of the {\greektext{p}} digit series with sequences of real non-uniform independent variables.}\label{sec:RandReSeq}

The top of Table \ref{tab:GausPoiExp} presents data on a string of random independent variables distributed following a Gauss $ f(x) = N(0,1) $ pdf. As seen, all values of $ D $ differ statistically from the chain of $ \pi $ digits when $ N'\geqq100 $ long. For a Gaussian variable outliers such as $\abs{x} \rightarrow \infty $ becomes likelier as $N' \rightarrow \infty$ which slows the convergence towards $ D=2 $ of a Gaussian sequence embedded in a unit square.

The middle of Table \ref{tab:GausPoiExp} presents data on a sequence of random independent variables exponentially distributed as
\begin{equation}\label{eq:Expon}
f_e(x|\eta)= \eta e^{-x\eta}.
\end{equation}
The exponential variables were generated as indicated in Section \ref{sec:Random},  which determines that $\eta=1$. For variables distributed as Eq. (\ref{eq:Expon}), $ \eta=1 \implies \mu= \sigma^2 = 1 $. As seen, sequences had $ D $ statistically distinct from the series of $\pi$ digits when $ N' \geqq  100$. As in the case of the Gaussian data, extreme values $ x \rightarrow \infty $ become likelier as $N' \rightarrow \infty$ which slows the convergence towards $ D=2 $ for this sequence.

The bottom of Table \ref{tab:GausPoiExp} presents data on a sequence of random independent variables distributed as an exponential pdf of the form
\begin{equation}\label{eq:Pois}
f_p(x|\eta)= \dfrac{\eta^x}{x!} e^{-\eta}.
\end{equation}
When simulating Poissonian variables as in Section \ref{sec:Random}, implicitly $ \eta = 1 \implies \mu= \sigma^2 = 1 $. From data in the table random Poissonian $ f(x|\eta=1) $ sequences had $ D $ which differs from $ \pi $ digit succession when $ N'\geqq  100,000$.

\begin{table}[h!]
\begin{center}
\caption{{Comparing $ D $ for $ \pi $ digits chains as $ N' $ grows from 10 to $ 10^{9} $ with $ D $ calculated for three different set of 30 periodic sequences.}}
\begin{tabular}{*{1}{m{6cm}}*{1}{m{3.5cm}}*{1}{m{1cm}}*{1}{m{1.3cm}}c}
\hhline{*{5}{=}}
\centering $ N'$ &\centering $ D \pm sd $ & \centering $ \lambda$ & \centering $ \epsilon $& $ \varPsi(\lambda) $\\
\hhline{*{5}{-}}
\centering Period of 300/397 \textit{n}=99&&&&\\
\centering 10&\centering 1.175190 \textpm{} 0.114982&0.474 & \centering NA&0.318\\
\centering 100&\centering 1.553983 \textpm{} 0.016619&0.258 & \centering NA&0.398 \\
\centering 1,000&\centering 1.688931 \textpm{} 0.003614&0.237 & \centering NA&0.406 \\
\centering 10,000&\centering 1.761546 \textpm{} 0.000875&2.919 & \centering 0.052 &$ 2 \cdot 10^{-3} $\\
\centering 100,000&\centering 1.806526 \textpm{} 0.000225&10.16 &\centering $ 4 \cdot10^{-3} $&$ <10^{-14} $\\
\centering 1,000,000&\centering 1.837233 \textpm{} 0.000060 &31.82 &\centering $ 4 \cdot 10^{-4} $&$ \ll10^{-14} $\\
\centering 10,000,000&\centering 1.859527 \textpm{} 0.000016 &102.7 &\centering $ 4 \cdot 10^{-5} $&$ \ll10^{-14} $\\
\centering 100,000,000&\centering 1.8764493 \textpm{} 0.0000045 &324.3 &\centering $ 4 \cdot 10^{-6} $&$ \ll10^{-14} $\\
\centering 1,000,000,000&\centering 1.8897328 \textpm{} 0.0000013 &1025 &\centering $ 4 \cdot 10^{-7} $&$ \ll10^{-14} $\\
\centering Period of 991/997 \textit{n}=166&&&\\
\centering 10&\centering 1.305775 \textpm{} 0.088497&0.859 & \centering NA&0.195 \\
\centering 100&\centering 1.543217 \textpm{} 0.017135&0.879 & \centering NA&0.190\\
\centering 1,000&\centering 1.678382 \textpm{} 0.003853&2.960 & \centering 0.051 &$ 2 \cdot 10^{-3} $\\
\centering 10,000&\centering 1.753121 \textpm{} 0.000933&6.291 & \centering 0.011&$ 2 \cdot 10^{-10} $\\
\centering 100,000&\centering 1.799649 \textpm{} 0.000240&19.19 &\centering $ 10^{-3} $&$ <10^{-14} $\\
\centering 1,000,000&\centering 1.831443 \textpm{} 0.000064 &61.03 &\centering $ 10^{-4} $&$ \ll10^{-14} $\\
\centering 10,000,000&\centering 1.854530 \textpm{} 0.000017 &191.1 &\centering $ 10^{-5} $&$ \ll10^{-14}$\\
\centering 100,000,000&\centering 1.8720543 \textpm{} 0.0000048 &604.6 &\centering $ 10^{-6} $&$ \ll10^{-14}$\\
\centering 1,000,000,000&\centering 1.8858104 \textpm{} 0.0000014 &1913 &\centering $ 10^{-7} $&$ \ll10^{-14}$\\
\centering Period of 1,000,001/999,997 \textit{n}=1508 &&&\\
\centering 10&\centering 1.105105 \textpm{} 0.121116&1.029 & \centering NA&0.152\\
\centering 100&\centering 1.522014 \textpm{} 0.020681&1.753 & \centering 0.145 &0.040\\
\centering 1,000&\centering 1.683780 \textpm{} 0.003751&1.601 & \centering NA &0.055\\
\centering 10,000&\centering 1.756380 \textpm{} 0.000925&2.822 & \centering  0.056&$ 2 \cdot 10^{-3} $ \\
\centering 100,000&\centering 1.802223 \textpm{} 0.000238&8.501 &\centering $ 6 \cdot 10^{-3} $&$ <10^{-14} $\\
\centering 1,000,000&\centering 1.833610 \textpm{} 0.000063 &27.22 &\centering $ 6 \cdot 10^{-4} $ &$ \ll10^{-14} $\\
\centering 10,000,000&\centering 1.856400 \textpm{} 0.000017 &84.09 &\centering $ 6 \cdot 10^{-5} $&$ \ll10^{-14} $ \\
\centering 100,000,000&\centering 1.8736994 \textpm{} 0.0000048 &266.3&\centering $ 6 \cdot 10^{-6} $ &$ \ll10^{-14} $\\
\centering 1,000,000,000&\centering 1.8872786 \textpm{} 0.0000016 &842.8&\centering $ 6 \cdot 10^{-7} $ &$ \ll10^{-14} $\\

\hhline{*{5}{=}}
\end{tabular}\label{tab:PeriodvsPi}
\end{center}
\footnotesize{{$D$: Calculated with Eq. \ref{eq:DCS} ; sd: standard deviation of $D$ calculated with Eq. \ref{eq:VarDemp} ; $N' $: number of $ \pi $ digits used to calculate $D$; $ \mathbb{Z}[0,9] $calculated as indicated in Section \ref{sec:Random} ; $ \lambda $: Calculated as in Eqs. \ref{eq:VP} and \ref{eq:PDif} for the difference between $D$ before and after of randomizing $ \pi $; $ \epsilon $: Probability calculated with theVysochanskij-Petunin inequality (Eq. \ref{eq:PVP}) for the difference in $D$ before and after of randomizing {\selectlanguage{greek}p}; NA: Not applicable. Values for $\pi$ in Table \ref{tab:PiPiRand}.}}
\end{table}

\subsection{\textit{D} properties of some rational number period digit sequences.}\label{sec:PeriodRation}

Table \ref{tab:PeriodvsPi} presents data on $D$ digit strings for three arbitrarily chosen periodic rational numbers. The numbers were:
\begin{equation}\label{eq:Periods}
\begin{split}
\text{Top of the table:}\quad \dfrac{300}{394}\;\quad\quad\quad\quad n= 99 \\
\quad\text{Middle of the table:}\quad \dfrac{991}{997}\;\;\;\quad\quad\quad n= 166 \\
\text{Bottom of the table:}\;\dfrac{1,000,001}{999.997}\quad n= 1508
\end{split}
\end{equation}
where $n$ is the number of decimals in the period. Decimal period was determined with the RealDigits[] function \citep{Wolfram2003} in \textit{Wolfram} Mathematica\textsuperscript{\textregistered} (Wolfram Research Inc.), and its digits were repeated as necessary to build chains of $10, 10^2, \ldots , 10^8, 10^9$ elements. These sequences were compared with the corresponding not randomized $\pi$ digit string as shown in Table \ref{tab:PiPiRand}. All periodic digit sequences in Table \ref{tab:PeriodvsPi}  were highly significantly different in $ D $ from non-randomized $ \pi $ digit sequence calculated when $N' \geqq 1000 $, at $ P(\lambda)=\epsilon < 0.05 $; if data would have been Gasussian [indicated in the table as $ \varPsi(\lambda) $], the significance would have been higher and the decision level ($ \alpha\leqq 0.05 $) would have been achieved with smaller $ N' $s.

\begin{table}[h!]
\begin{center}
\caption{{Comparing $ D $ for  $ \pi $ digits chain as $ N' $ grows from 10 to $ 10^{9} $ with $ D $ calculated for three different set of 30 periodic sequences after randomizing them.}}
\begin{tabular}{*{1}{m{8cm}}*{1}{m{3.5cm}}*{2}{m{1.4cm}}c}
\hhline{*{5}{=}}
\centering $ N'$ &\centering $ D \pm sd $ & \centering $ \lambda$ & \centering $ \epsilon $& $ \varPsi(\lambda) $\\
\hhline{*{5}{-}}

\centering Randomized period of 300/397 \textit{n}=99&&&&\\
\centering 10 &\centering 1.295488 \textpm{} 0.087319&\centering  8.176& \centering $6 \cdot 10^{-3}$ & $  <10^{-14} $\\
\centering 100&\centering 1.524229 \textpm{} 0.017104 &\centering  63.22& \centering $10^{-3}$  &$ \ll10^{-14} $ \\
\centering 1,000&\centering 1.687801 \textpm{} 0.003580&\centering 331.1& \centering $4 \cdot 10^{-6}$&$ \ll10^{-14} $\\
\centering 10,000&\centering 1.760949 \textpm{} 0.000882&\centering 1417& \centering $2 \cdot 10^{-7}$ &$ \ll10^{-14} $\\
\centering 100,000&\centering 1.806274 \textpm{} 0.000225&\centering 5678&\centering $ 10^{-8} $&$ \ll10^{-14} $\\
\centering 1,000,000&\centering 1.836768 \textpm{} 0.000060 &\centering $ \approx 2 \cdot 10^{4}  $&\centering $ 8 \cdot 10^{-10} $&$ \ll10^{-14} $\\
\centering 10,000,000&\centering 1.859150 \textpm{} 0.000016&\centering  $ \approx 8 \cdot 10^{5}  $&\centering $ 7\cdot 10^{-11} $&$ \ll10^{-14} $\\
\centering 100,000,000&\centering 1.8761116 \textpm{} 0.0000045&$ \approx 3 \cdot 10^{5}  $ &\centering $ 5 \cdot 10^{-12} $&$ \ll10^{-14} $\\
\centering 1,000,000,000&\centering 1.8894330 \textpm{} 0.0000013 &\centering  $ \approx 10^{6}$ &\centering $ 4 \cdot 10^{-13} $&$ \ll10^{-14} $\\
\centering Randomized period of 991/997 \textit{n}=166&&&&\\
\centering 10 &\centering 1.333458 \textpm{} 0.08690 & \centering10.04 & \centering $4 \cdot 10^{-3}$ & $  <10^{-14} $\\
\centering 100&\centering 1.568944 \textpm{} 0.016129 & \centering 65.94 & \centering $10^{-4}$  &$ \ll10^{-14} $ \\
\centering 1,000 &\centering 1.681535 \textpm{} 0.003727& \centering 313.0& \centering  $5 \cdot 10^{-6}$ &$ \ll10^{-14} $\\
\centering 10,000&\centering 1.759015 \textpm{} 0.000882& \centering 1369& \centering $2 \cdot 10^{-7}$&$ \ll10^{-14} $\\
\centering 100,000&\centering 1.804497 \textpm{} 0.000230& \centering 5483 &\centering $ 10^{-8} $&$ \ll10^{-14} $\\
\centering 1,000,000&\centering 1.835480 \textpm{} 0.000060& \centering $ \approx 2 \cdot 10^{4} $ &\centering $ 10^{-11} $&$ \ll10^{-14} $\\
\centering 10,000,000&\centering 1.858072 \textpm{} 0.000016 & \centering $ \approx 8 \cdot 10^{4} $&\centering $ 7 \cdot 10^{-11} $&$ \ll10^{-14}$\\
\centering 100,000,000&\centering 1.8751709 \textpm{} 0.0000046& \centering $ \approx 3 \cdot 10^{5}$&\centering $ 6 \cdot 10^{-12} $&$ \ll10^{-14}$\\
\centering 1,000,000,000&\centering 1.8885890 \textpm{} 0.0000013& \centering $ \approx 10^{6}$ &\centering $ 4 \cdot 10^{-13} $&$ \ll10^{-14}$\\
\centering Randomized period of 1000001/999997 \textit{n}=1508 &&&&\\
\centering 10&\centering 1.111450 \textpm{} 0.165965&\centering 4.820& \centering 0.019 &$ 7 \cdot 10^{-7} $\\
\centering 100&\centering 1.551871 \textpm{} 0.01747&\centering 56.56& \centering $10^{-4} $&$ <10^{-14}$\\
\centering 1,000&\centering 1.691323 \textpm{} 0.003627&\centering 323.4& \centering $ 4 \cdot 10^{-6} $ &$ \ll10^{-14}$\\
\centering 10,000&\centering 1.761373 \textpm{} 0.000888 &\centering 1373& \centering  $ 2 \cdot 10^{-7} $&$ \ll10^{-14}$ \\
\centering 100,000&\centering 1.806941 \textpm{} 0.000225&\centering 5517 &\centering $ 10^{-8} $&$ \ll10^{-14} $\\
\centering 1,000,000&\centering 1.837313 \textpm{} 0.000060 &\centering $ \approx2 \cdot 10^{4} $ &\centering $ 10^{-9} $ &$ \ll10^{-14} $\\
\centering 10,000,000&\centering 1.859554 \textpm{} 0.000016 &\centering $\approx 8 \cdot 10^{4} $&\centering $ 7 \cdot 10^{-11} $&$ \ll10^{-14} $ \\
\centering 100,000,000&\centering 1.8764907 \textpm{} 0.0000046&\centering $\approx 3 \cdot 10^{5} $&\centering $ 6 \cdot 10^{-12} $ &$ \ll10^{-14} $\\
\centering 1,000,000,000&\centering 1.8897681 \textpm{} 0.0000013 &\centering $\approx 10^{6} $&\centering $ 4 \cdot 10^{-13} $ &$ \ll10^{-14} $\\
\hhline{*{5}{=}}
\end{tabular}\label{tab:RandPervsPer}
\end{center}
\footnotesize{{$D$: Calculated with Eq. \ref{eq:DCS} ; sd: standard deviation of $D$ calculated with Eq. \ref{eq:VarDemp} ; $N' $: number digits used to calculate $D$; $ U_{\mathbb{Z}}[0,9] $calculated as indicated in Section \ref{sec:Random} ; $ \lambda $: Calculated as in Eqs. \ref{eq:VP} and \ref{eq:PDif} for the difference between $D$ before and after of randomizing $ \pi $; $ \epsilon $: Probability calculated with the Vysochanskij-Petunin inequality (Eq. \ref{eq:PVP}) for the difference in $D$ before and after of randomizing {\selectlanguage{greek}p}; NA: Not applicable. Periodic rational number digit chains were randomized as indicated in Section \ref{sec:PeriodRation}.$\pi \; D$ values in Table \ref{tab:PiPiRand}.}}
\end{table}

To study the effect of randomization on the periodic sequences described in the previous paragraph, the digits in the period of length $n$ were randomly sampled as $ U_{\mathbb{Z}}[1,n] $ and the digit obtained in each sampling was used to build sequence initially of 10 digits, which was expanded by adding 90 digits to this sequence to make one of 100 digits, and so on up to $ N' =10^9 $. Thirty sets of randomized series were calculated for each one of the periodic digit sequences in equation set (\ref{eq:Periods}), $D$ was calculated for all of those 30 set, $N'$ started at 10 and was incremented in steps of $\times 10$ up to $10^9$, as before. $ \overline{D} \pm sd $ of the randomized periodical digit catenae is presented in Table \ref{tab:RandPervsPer} together with the values of $ \lambda $, $ \varepsilon $ and $ \varPsi(\lambda) $ obtained when $ \overline{D} \pm sd  $ of 30 sets of integers were compared with the corresponding native unhashed sequences in Table \ref{tab:PeriodvsPi}{ }.

Table \ref{tab:RandPervsPer} shows that randomization modified the $D$ values of the periodic sequence in highly statistically significant manner; the changes and their significances were similar after randomizing decimals of the three periodic rational numbers. In all cases the change was an increase in $D$ as it would be expected if randomization induced disorder and increased their entropy (More details in Section \ref{sec:Intro}).
To get an insight on the characteristics of the digit sequences of the rational number periods in equation set (\ref{eq:Periods}), the statistical properties of their digit sets we studied. For this purpose a set of digits of one period from the three rational number,s were considered as set of consecutive random variables, and statistically compared with similar sets of the other two periodic fractions. This analysis  indicated that they were all distributed as $ U[0,9] $. The nonparametric Kruskall-Wallis analysis of variance failed to detect any difference between them ($ P(\text{H0}) = 0.490 $). This was confirmed by the Mann-Whitney (Wilcoxon) test ($0.95 \geqq P(\text{H0}) > 0.85$). The Smirnov test \citep{Smirnov1939, Smirnov1939a} based on Kolmogorov statistics  \citep{Kolmogorov1933} was unable to find differences between the distributions of digits in the three periods ($ 1 \geqq P(\text{H0})  > 0.61$). 

All statistical test mentioned in the previous paragraph consider the digits as unordered sets. To test if there is some trend within the digits of a rational number unshuffled period, the set was studied with the Wilks \citep{Guttman1971} above and below the median test; this test is an analysis of runs, sequences of consecutive identical events above or bellow the median value; for the purpose of the test the minimum run length is 1 for a single isolated event above or below the median. The Willks tests indicated that period\textquoteright{s} digit sequences were randomly distributed about their median value [$ 0.78 > P(\text{H0})  > 0.38 $]: i.e., they had no trend. All these test indicate that the digits of the periods have the same statistical properties, as well as that they not innerly ordered in any way, which suggests that the rational numbers considered, may be normal numbers.

\begin{table}[h!]
\begin{center}
\caption{Values of $ D $ for two random walks.}
\begin{tabular}{*{1}{m{5cm}}*{1}{m{3.5cm}}*{1}{m{1.6cm}}cr}
\hhline{*{4}{=}}
\centering $ N'$ &\centering $ D \pm sd $ & \centering $ \lambda$ & \centering $ \epsilon(\lambda) $&\\
\hhline{*{4}{-}}

\centering Brownian $N(0,1)$ random walk&&&\\
\centering 10 &\centering 1.096039 \textpm{} 0.098634 &&\\
\centering 100 &\centering 1.267443 \textpm{} 0.056499&&\\
\centering 1,000 &\centering 1.328054 \textpm{} 0.035347&&\\
\centering 10,000 &\centering 1.371589 \textpm{} 0.029656&&\\
\centering 100,000 &\centering 1.384101 \textpm{} 0.023168&&\\
\centering 1,000,000 &\centering 1.400816 \textpm{} 0.016508&&\\
\centering 10,000,000 &\centering 1.420885 \textpm{} 0.020553&&\\
\centering 100,000,000 &\centering 1.427627 \textpm{} 0.013448&&\\
\centering 1,000,000,000 &\centering 1.435309 \textpm{} 0.010149&&\\
\centering Uniform $U[-1,1]$ random walk&&&\\
\centering 10 &\centering 1.087430 \textpm{} 0.132544&\centering 0.052&NA\\
\centering 100 &\centering 1.236725 \textpm{} 0.054427&\centering 0.391&NA\\
\centering 1,000 &\centering 1.319483 \textpm{} 0.035465&\centering 0.171&NA\\
\centering 10,000 &\centering 1.364712 \textpm{} 0.032263&\centering 0.157&NA\\
\centering 100,000 &\centering 1.391049 \textpm{} 0.024396&\centering 0.207&NA\\
\centering 1,000,000 &\centering 1.403320 \textpm{} 0.017006&\centering 0.106&NA\\
\centering 10,000,000 &\centering 1.421456 \textpm{} 0.018172&\centering 0.021&NA\\
\centering 100,000,000 &\centering 1.429921 \textpm{} 0.014729&\centering 0.115&NA\\
\centering 1,000,000,000 &\centering 1.436430 \textpm{} 0.013981&\centering 0.065&NA\\

\hhline{*{4}{=}}
\end{tabular}\label{tab:Brown}
\end{center}
\footnotesize{{$D$: Calculated with Eq. \ref{eq:DCS} ; sd: standard deviation of $D$ calculated with Eq. \ref{eq:VarDemp} ; $N' $: number digits used to calculate $D$. The two random walks were $\{ r_{i +1} = r_i + N(0,1)\}$ and $\{ r_{i+1} = r_i + U[-1,1)\}$. The values of $ \lambda $ and $ \epsilon(\lambda) $ in this table, correspond to comparisons of significance of the two sequences in this table, See section \ref{sec:Random} or other details.}}
\end{table}

\subsection{On some sequences of non-independent events.}

Table \ref{tab:Brown} presents data on a different kind of sequences, they are usually known as \textit{random walks}, sets of dependent random variables which will be denoted as $ \{r_i\} $. They can be described as recursions where $\{ r_{i+1} \looparrowleft r_i + g(\{\theta_k\} \}_{i=1,2,\ldots,N-1}$ where $ g(\{\theta_k\}) $ is a random process dependent on a set of parameters $\{\theta_k\}_{k=1,2, \ldots, m} $. Two such functions [See Eqs, (\ref{eq:BrownW}) and (\ref{eq:UnifW})] are presented in Table \ref{tab:Brown} .

Three characteristics are evident frmm $ \overline{D} $ values in Table \ref{tab:Brown}, First: $ D $ seems to converge towards $ D\leqq 1.5 $ which is the fractal dimension of Brownian noise \citep{Mandelbrot1975, Mandelbrot1983}. Second: standard deviations of $ D $ in Table \ref{tab:Brown} do decrease as $ N' $ increases, but the decrease is much slower than for any of de sequences of independent random variables (and of $ \pi $ too) previously considered here. Third: although the $ D $ values are obviously different from any of the independent sequences of random variables discussed previously in this work, the corresponding (same $ N' $) value of the Brownian and uniform random walks in the Table \ref{tab:Brown} seem identical, at least within the range of $ N' \leqq 10^9$ studied.

\section{Discussion.}

\subsection{On waveforms and digit sequences considered as waveform samples.}

The term waveform applies to planar curves shaped as a wave, usually drawn as instantaneous values of a changing quantity ${\textit{versus}}$ time. Any waveform is an infinite series of points. Aside of classical methods such as moment statistics and regression analysis, properties such as the Kolmogorov-Sinai-Chaitin entropy \citep{Sinai1959, Kolmogorov1965, Chaitin1969, Grassberger1983}, the apparent entropy \citep{Pincus1991} and the fractal dimension \citep{Sevcik1998a} have been proposed to tackle the problem of pattern analysis of waveforms. 

When analyzing waveforms it is usually convenient to consider finite samples of points $ (x_i, y_i) $ separated constant increments $ \varDelta x = x_{i+1} -x_i $, which are usually transformed (or \textit{discretized}) into a set $ \{z_i \}$ where $ y_i \mapsto z_i  \in \mathbb{Z}$  for which $ 2^{-\delta} \leqq z_i \leqq 2^\delta $ where $ \delta \in \mathbb{Z} $ is a constant determining the \textit{resolution} of the discretization \citep{Smith1997}.

A simple algorithm exists (see Eq. \ref{eq:DCS}) able to approximate the fractal dimension of a discretized waveform; and if $N'\rightarrow \infty$, $D \rightarrow \Phi$ the curve\textquoteright{s} true fractal dimension \ref{eq:DCS}. Furthermore, \citep{Sevcik1998a}, an interesting finding is also that Eq. (\ref{eq:DCS}) applied to sets of \textit{N{\textquotesingle}} samples may distinguish between waveforms that seemingly converge to the same value when ${N' \rightarrow \infty}$. This is the case of the uniform and Gaussian white noises. If  $ u $ and $ v $ are \textit{U[0,1]} random variables, then ${\underset{N' \rightarrow \infty }\lim  f(u - v) \approx{N \left( \textit{0, 1/6} \right) }}$ and the expectation $\nu=\mathbb{E}\abs{u-v }\approx \textit{1/3}$:
\begin{equation}\label{eq:Dranseq}
D= 1+\underset{N' \rightarrow \infty }{\lim } \left \{ \frac{ \ln \left[N'\cdot \sqrt{\nu ^{2}+\frac{1}{{N'}^2}} \right] }{ \ln (2 \cdot N')} \right \} = 2
\end{equation}
which holds also for any sequence $ \{y_i\}_{i=1,2,\ldots,N'}$ fulfilling 

\begin{equation}\label{eq:Lim2}
(y_i  \upmodels y_{k \neq i}) \forall y_{k \neq i}  \implies {\underset{N' \rightarrow \infty} \lim} D = 2 
\end{equation}

\noindent{alternatively $ D $ will also converge towards a constant when $ N' \rightarrow \infty $} if
\begin{equation}\label{eq:LimNot2}
(y_i  \nupmodels y_{k \neq i}) \forall y_{k \neq i}  \implies {\underset{N' \rightarrow \infty} \lim} D < 2
\end{equation}
but the way $ D $ reaches a limit as $ N' \rightarrow \infty $ depends on the properties of the  sequence, as it is shown by data in Tables \ref{tab:PiPiRand} through \ref{tab:Brown}. Under both conditions expressed by Eqs. (\ref{eq:Lim2}) and (\ref{eq:LimNot2}) $\underset{N' \rightarrow \infty}{\lim}{\Var(D)}= 0$ (see Eq. \ref{eq:ConvVarD}). 

The sequence of $ \pi $ digits may be visualized as a low resolution sample $ \{0 \leqq z_i \leqq 9  \}$ from a hypothetical waveform. This waveform was linearly transformed by embedding the $ \pi $ digit string it into a unit square as indicated in Eqs. (\ref{eq:xtrans}) and (\ref{eq:ytrans}). In the unit square the digits were connected by $ N' $ straight line segments, mimicking a discretization of the hypothetical discretized waveform, and the length of these segments was calculated. The hypothetical waveform\textquoteright{s} fractal dimension was approximated with Eq. (\ref{eq:DCS}) \citep{Sevcik1998a} using $ N' $ values, starting at 10, and increasing in steps of 10 up to $ 10^9 $.

\subsection{Sequence randomization and fractal dimension.}\label{sec:BoltzEnt}
Randomness means disorganization, if a system is truly random it cannot be \textquotedblleft{}more random\textquotedblright{}. In Boltzmann terms \citep{Boltzmann1964}, the entropy $ S $ and information ($ \iota $) are defined as 
\begin{equation}\label{eq:Boltzmann}
\begin{split}
S = k_B \ln (W) \quad  \Leftrightarrow  \quad \iota = -k_B \ln (W).
\end{split}
\end{equation}
Where $k_B$ is Boltzmann\textquoteright{}s constant and $ W $ is the number of states in a system and $ \iota $ is negentropy \citep{Szilard1929, Shannon1948}, negative entropy. As described in Section \ref{sec:Intro}, the most general interpretation of entropy is as a measure of our uncertainty about a system, in a random sequence of size $ \{u_i\} \supseteq \varOmega$ has a maximum entropy, any $ \{ u_i\} $ cannot be further disorganized. For a sequence or a waveform, as implied by Ec. \ref{eq:Dranseq} ,
\begin{equation}\label{eq:S_D}
\underset{S \rightarrow \infty }{\lim }D=2.
\end{equation}

Data in Table \ref{tab:PiPiRand} shows that $ D $ did not change when $ \pi $ sequence was randomized. The upper half of the table presents $ \pi\; D$ values calculated with Eq. \ref{eq:DCS} and their standard deviation calculated with Eq. \ref{eq:VarD}. The lower part of the table presents the mean $ D \pm s.d.$ calculated for 30 series of $ \pi $ digits randomized as indicated in Section \ref{sec:Random}. $D \rightarrow 2$ and $ \Var(D)\rightarrow 0 $ as $N' \rightarrow\infty $ in either native or randomized $ \pi $ digit sequences. The most striking finding is that when comparing $ D $ calculated with the same $ N' $ the values are equal before and after randomization.

This is to be expected if $ \pi $ digit sequence is random and uniform, if the probabilities of all digits $ p(0)=p(1)=\ldots=p(9)=0.1 $, and Eq. (\ref{eq:Lim2}) holds; i.e. the sequence is not deterministic in any manner. Given any digit in the sequence, it provides no information ($\iota$) on what happened before and after in the sequence. \emph{Those are the properties of a normal number}.

Data in Table \ref{tab:PiVsU09} is similar to Table \ref{tab:PiPiRand}, but it shows data for $ D $ of 30 series of  $ U_{\mathbb{Z}}[0,9] $ digits (Upper part of the table) are not different from the $ \pi $ digit series while data for $ D $ of 30 series of   $ U[0,1] $ real numbers (Lower part of the table) are highly distinct from $ \pi $ values. A common feature of both types of chains is that $D \rightarrow 2$ and $ \Var(D)\rightarrow 0 $ as   $N' \rightarrow\infty $. These results show the ability of Eq. \ref{eq:DCS} to differentiate between distinct series of independent uniform random numbers, and reinforces the notion that $ \pi $ digit sequence corresponds to a normal number.

\subsection{Sequence of {\greektext{p}} compared with independent Gaussian, exponential and Poissonian random variables\textquoteright{} sequences.}

Table \ref{tab:GausPoiExp} presents $ D $ data for another three series of real number sequences distributed a a Gaussian $ N(0,1) $ (top of the table),  exponential $ f_e(x|\eta=1)$ (middle of the table) and Poissonian   $f_p(x | \eta=1) $ (bottom of the table) independent random real numbers sequences (see \ref{sec:Random} and \ref{sec:RandReSeq} for details on random variable generation and the specific pdf functions). As shown all three types of series are different from the $ \pi $ digits sequence although for all of them $D \rightarrow 2$ and $ \Var(D)\rightarrow 0 $ as   $N' \rightarrow\infty $.

\subsection{Comparing random walks}

Table \ref{tab:Brown} presents data to compare two random walks, one is the Brownian walk  where $\{ r_{i+1} = r_i + N(0,1)$  and the other, called Uniform $ U[-1,1] $, where $\{ r_{i+1} = r_i + U[-1,1]\}$. As seen in the table this two walks are indistinguishable when $ N' \leqq 10^9 $, and converge to $ D \leqq 1.5 $. The fractal dimension of Brownian noise is 1.5 \citep{Mandelbrot1975, Mandelbrot1983}. The data coincides with previous results \citep{Sevcik1998a} which indicate that fractal analysis of time series using Ec. (\ref{eq:DCS}) converges toward the right $ \Phi $ value when random walks are studied. The data shows that $ \pi $ digits are indeed not a random walk of integers. 

\subsection{Fractal analysis of periodic sequences.}\label{sec:StrucDigPerSec}

Equation \ref{eq:DCS} was used to calculate $ D $ of digit sequences which contains some structure. For this purpose, sequences of up to $ 10^9 $ decimal digits of three rational numbers with periods between 99 and 1508 digits were considered (see Section \ref{sec:PeriodRation} for details). As shown in Table \ref{tab:PeriodvsPi} the sequences of $ D $ values obtained with the 3 periodic series were clearly distinct from the $ \pi $ succession of digits. This shows that periodic structure in decimal digits may be detected and suggests, once again, no structure exists in the decimal digits of $ \pi $,  as it should be if $ \pi $ is a normal irrational number.

\begin{table}[h!]
\begin{center}
\caption{{Statistical significance of differences between  $ D $ values of the series of uniform normal digits of type $ U_{\mathbb{Z}}[0,9] $ and three periodic digit sequUences digits as $ N' $ grows from 10 to $ 10^{9} $.}}
\begin{tabular}{*{1}{m{5.5cm}}*{1}{m{3cm}}*{2}{m{4cm}}c}
\hhline{*{5}{=}}
&\centering $U_{\mathbb{Z}}[0,9]$ &  \centering 1,000,001/999,997 period& \centering 991/997 period &\\
\centering $ N'$ &\centering $ \epsilon(\lambda) $ & \centering $\epsilon(\lambda)$ & \centering $ \epsilon(\lambda) $&\\
\hhline{*{5}{-}}

\centering 300/397 period \textit{n}=99&&&\\
\centering 10&\centering NA&\centering NA&\centering NA&\\
\centering 100&\centering 0.021&\centering NA&\centering  NA& \\
\centering 1,000&\centering $ 4 \cdot 10^{-3} $&\centering NA&\centering 0.056&\\
\centering 10,000&\centering $ 2 \cdot 10^{-3} $&\centering 0.014&\centering $ 5 \cdot 10^{-3} $&\\
\centering 100,000&\centering $ 9 \cdot 10^{-4} $&\centering $ 10^{-3} $&\centering $ 5 \cdot 10^{-4} $&\\
\centering 1,000,000&\centering $ 6 \cdot 10^{-4} $&\centering $ 10^{-4} $&\centering  $ 5 \cdot 10^{-5} $&\\
\centering 10,000,000&\centering $ 2 \cdot 10^{-4} $&\centering $ 10^{-5} $&\centering $ 5 \cdot 10^{-6} $&\\
\centering 100,000,000&\centering $ 2 \cdot 10^{-4} $&\centering $ 10^{-6} $&\centering $ 5 \cdot 10^{-7} $&\\
\centering 1,000,000,000&\centering $ 10^{-4} $&\centering $ 10^{-7} $&\centering $ 5 \cdot 10^{-8} $&\\
\centering 991/997 period \textit{n}=166&&&\\
\centering 10&\centering NA&\centering NA&&\\
\centering 100&\centering 0.025&\centering NA&& \\
\centering 1,000&\centering $ 4 \cdot 10^{-3} $&\centering NA&&\\
\centering 10,000&\centering $ 2 \cdot 10^{-3} $&\centering 0.036&&\\
\centering 100,000&\centering $ 10^{-3} $&\centering $ 4 \cdot 10^{-3} $&&\\
\centering 1,000,000&\centering $ 7 \cdot 10^{-4} $&\centering $ 4 \cdot 10^{-4} $&&\\
\centering 10,000,000&\centering $ 3 \cdot 10^{-4} $&\centering $ 4 \cdot 10^{-5} $&&\\
\centering 100,000,000&\centering $ 2 \cdot 10^{-4} $&\centering $ 4 \cdot 10^{-6} $&&\\
\centering 1,000,000,000&\centering $ 10^{-4} $&\centering $ 4 \cdot 10^{-7} $&&\\
\centering 1,000,001/999,997 period \textit{n}=1508 &&&\\
\centering 10&\centering NA&&&\\
\centering 100&\centering NA&&& \\
\centering 1,000&\centering NA&&&\\
\centering 10,000&\centering 0.014&\centering &\centering &\\
\centering 100,000&\centering $ 10^{-3} $&\centering &\centering &\\
\centering 1,000,000&\centering $ 2 \cdot 10^{-4} $&\centering &\centering &\\
\centering 10,000,000&\centering $ 10^{-5} $&\centering &\centering &\\
\centering 100,000,000&\centering $ 10^{-6} $&\centering &\centering &\\
\centering 1,000,000,000&\centering $ 10^{-7} $&\centering &\centering &\\
\hhline{*{5}{=}}
\end{tabular}\label{tab:UniPervsPeri}
\end{center}
\footnotesize{$N' $: number of $ \pi $ digits used to calculate $D$; $ U_{\mathbb{Z}}[0,9] $ calculated as indicated in Section \ref{sec:Random}; $ \epsilon $: Probability calculated with theVysochanskij-Petunin inequality (Eq. \ref{eq:PVP}) for the difference in $D$ beween the different sequences tabulated; $ n $ is the rational umber period length; NA: Not applicable. Other details as in previous tables}
\end{table}

In the previous paragraph, three structured digit periodic sequences were with $ \pi $ digit sequence and were found to be different. In Table \ref{tab:RandPervsPer}, 30 sets of each of the periodic sequences considered in Table \ref{tab:PeriodvsPi} were generated by randomization as indicated in Section \ref{sec:PeriodRation} and compared with the native, non-randomized, sequence. In all three cases randomization increased $ D \rightarrow 2 $ in a highly statistically significant manner for sequences having $ N' \geqq 10 $, and the statistical significance was evident using the distribution independent Vysochanskij-Petunin Inequality. (Section \ref{sec:VP}). The finding shows that when structure is destroyed $ D$ get closer to 2, as predicted by Eq. \ref{eq:S_D}.

Table \ref{tab:UniPervsPeri} contains information on the statistical significance of differences between $ D $ values calculated for $ U_{\mathbb{Z}}[0,9] $ sequences and Brownian $ N(0,1) $ or uniform random walks with sequences of diverse $ N' $ values. Only values of $ \epsilon(\lambda) $ [Eq. (\ref{eq:PVP})] are presented. It is clearly seen in the table that Eq. \ref{eq:DCS} differentiates between the four types of decimal integer sequences when $ N' > n$ (where $ n $ is the period length), this suggests strongly that the periods of the rational numbers $ 300/397 $ ($ n=99 $), $ 991/997 $ ($ n=166 $) and 1,000,001/999,997 ($n=1508 $) are uniform random variables indistinct from $ U_{\mathbb{Z}}[0,9] $. Data in Table \ref{tab:UniPervsPeri} shows that  $ D $ calculated with Eq. (\ref{eq:PVP}) differentiates between sequences with different types of internal structure from a natural number with digits ordered as a $ U_{\mathbb{Z}}[0,9] $ chain.

\subsection{Is {\greektext{p}} normal?.}

The question has been asked by many authors \cite{Bailey2001, Bailey2012a}. The definition of normal number \citep{Borel1909} is \cite[pg. 299]{Sierpinski1988}:
\begin{quotation}
Let $ g $ be a natural number $ > 1 $; we write a real number $ x:x = [x] 
+(0.c1\, c2\,c3\ldots)_g$  as a decimal in the scale of $ g $. For any digit $ c $ (in the scale of $ g $) and every natural number $ n $ we denote by $ I (c, n) $ the number of those digits of the sequence $  c_1, c_2 , \ldots c_n$; which are equal to $ c $. If
\begin{equation*}
\underset{n}{\lim }\dfrac{I(c,n)}{n}=\dfrac{1}{g}
\end{equation*}
for each of the 9 possible values of $ c $, then the number $ x $ is called normal in
the scale of $ g $.
\end{quotation}
A number which is natural at any scale $ g $ is called \textit{absolutely natural} \citep{Becher2002}. For a number in base 10, the definition implies that $ c_i \in \mathbb{Z}[0,9] $ if it is to apply to any number as $ N \rightarrow \infty $ fractionary digits. Thus different authors appealed to statistical tests such as comparing the frequencies of each digit in the decimal sequence of $ \pi $, and the frequency of diverse combinations of the decimal digits \citep{Bailey2012a}. This approach is unsatisfactory since not all possible combinations can be evaluated, and is bound by the $ N $ value of the series studied, ans well as the size of the subsample of $ N $ used to do the statisticall tests.

The fractal analysis used in this work considered the series of $ \pi $ decimals, and calculated its approximate fractal dimension, for $ N= 10^1, 10^2,10^2,\ldots,10^{9} $ using Eq. (\ref{eq:DCS}) and showed that $ \underset{N \rightarrow \infty }{\lim} D \approx 2 $ for all series of type $ y_i = f(\{ \theta_k\}_i) $, where $ f(\{ \theta_k\}_i) $ is some random variable wich depends on a set of parameters. This was found true in this paper for series of real numbers distributed as pdf such as Gauss, Poisson, exponential or uniform ($ U[0,1] $) functions, as well as for discrete distribution such as $ U_{\mathbb{Z}}[0,9]$ and \emph{the decimal sequence of $ \pi $}. We have also seen that $ D $ for series obeying the condition represented by Eq. (\ref{eq:Lim2})  is constant under randomization, \emph{i.e.: does not change when the sequence is randomized and that this also a property of the digit sequence of $ \pi $}. 

Randomization increases the entropy (of the Boltzmann type (Section \ref{sec:BoltzEnt} and Eq. \ref{eq:Boltzmann}) \citep{Boltzmann1964}, or of the algorithmic type \citep{Sinai1959, Kolmogorov1965, Chaitin1969}) of a sequence, a random sequence has the maximum entropy or equivalently the minimum negentropy or information content \citep{Szilard1929, Shannon1948}. These conclusions could by falsified should singularity in the infinite series functions used to calculate $ \pi $ digits  \citep{Chudnovsky2000}, exist.
 
The value of $ D $ calculated for the $ \pi $ digits sequence was indistinguishable at all $ N' $ values used, from $ D $ values obtained for uniform random sequences of decimal digits ($ \mathbb{Z}[0,9] $), which are random by design (Section \ref{sec:Random}). Three periodic digits sequences were easily differentiated from $ \pi $\textquoteright{}s sequence and were not constant under randomization.
 
Infinity is far away. The findings in this study are based on the first $ 10^9 $ initial digits of $ \pi $; this limit was set by the computer power available to the author. One billion digits seems little given that 22,459,157,718,361 decimal and 18,651,926,753,033 hexadecimal digits \citep{Trueb2016a, Trueb2016b} of $ \pi $ have been calculated. Extending this study to the full known $ \pi $ digit sequence is only a problem of storage and computer power. Yet, the most interesting aspect of this study is that the more series digits you consider, the stronger its conclusions get. Unless the $ \pi $ digit sequence has a singularity and becomes not representable by continuous infinite series such as the Ramanujan series \citep{Berggren2000}, all the properties of decimal digits must be as shown here. Fractal analysis presented here indicates that $ \pi $ is normal in base 10. 
 
 \section{Acknowledgments.}

This manuscript was written in \LaTeX{ }using \textit{{\TeX}studio} for Linux (Also available for Apple OS X and MS Windows, http://www.texstudio.org), an open source free  \LaTeX{ }editor.


\end{document}